\documentclass[12pt]{amsart}
\usepackage{verbatim}
\usepackage[draft]{todonotes}

\usepackage{amsthm}
\usetikzlibrary{arrows}
\usepackage{blindtext}

\usepackage{amssymb}
\usepackage{mathrsfs}
\usepackage[all]{xy}
\usepackage{tikz}
\usepackage[alphabetic,nobysame,abbrev]{amsrefs}
\usepackage{enumerate}
\RequirePackage{fix-cm}
\usepackage{bm}
\usepackage{mdframed}
\usepackage{cancel}
\usepackage{makecell}
\usepackage{diagbox}

\usepackage{amscd}
\usepackage{amsmath}
\usepackage{bm}
\usepackage{newtxtext} 
\numberwithin{equation}{section}
\usepackage{romanbar}
\usepackage[toc,page]{appendix}
\usepackage{hyperref}
\input xy
\xyoption{all}
\usepackage{color}
\RequirePackage{fix-cm}
\usepackage{bm}
\usepackage{mdframed}
\usepackage{cancel}
\usepackage{makecell}
\usepackage{diagbox}
\usepackage{pgf,tikz}
\usepackage{mathrsfs}
\usetikzlibrary{arrows}
\usetikzlibrary{arrows,decorations.markings}
\usetikzlibrary{matrix}
\usetikzlibrary{graphs}
\usetikzlibrary{backgrounds}

\definecolor{qqqqff}{rgb}{0.,0.,1.}
\definecolor{xdxdff}{rgb}{0.49019607843137253,0.49019607843137253,1.}

\definecolor{qqqqff}{rgb}{0.,0.,1.}

\usepackage{tikz}
\usepackage{graphicx}
\usepackage{hyperref}
\usepackage[utf8]{inputenc}
\usepackage{amsmath,amsfonts,amssymb,euscript,graphicx,color,tikz}

\newtheorem{lemma}{Lemma}[section]
\newtheorem{corollary}[lemma]{Corollary}
\newtheorem{theorem}[lemma]{Theorem}
\newtheorem{proposition}[lemma]{Proposition}

\theoremstyle{definition}
\newtheorem{remark}[lemma]{Remark}
\newtheorem{notation}[lemma]{Notation}
\newtheorem{definition}[lemma]{Definition}

\newtheorem{example}[lemma]{Example}

\DeclareMathOperator{\Mod}{Mod}
\DeclareMathOperator{\modd}{mod}

\DeclareMathOperator{\Add}{Add}

\DeclareMathOperator{\Hom}{Hom}

\DeclareMathOperator{\Ab}{Ab}

\DeclareMathOperator{\Ker}{Ker}
\DeclareMathOperator{\Coker}{Coker}

\DeclareMathOperator{\Imm}{Im}

\DeclareMathOperator{\proj}{proj}

\newtheorem{question}[lemma]{Question}
\newtheorem*{theorem 0*}{Theorem}
\newtheorem*{theorem a*}{Theorem A}
\newtheorem*{theorem b*}{Theorem B}

\newcounter{diagram}
\numberwithin{diagram}{section}

\begin{document}
	
	\title{On definable subcategories}
	
	\author{Ramin Ebrahimi}
	\address{School of Mathematics, Institute for Research in Fundamental Sciences (IPM), P.O. Box: 19395-5746, Tehran, Iran}
	\email{ramin.ebrahimi1369@gmail.com / r.ebrahimi@ipm.ir}

	\subjclass[2020]{{16Gxx}, {16Dxx}, {16Yxx}, {18Dxx}, {18Bxx}}
	
	\keywords{Definable subcategories, The free abelian category, Auslander-Gruson-Jensen duality}

	\begin{abstract}
		Let $\mathcal{X}$ be a skeletally small additive category. 
		Using the canonical equivalence between two different presentations of the free abelian category over $\mathcal{X}$, we give a new and simple characterization of definable subcategories of $\Mod\text{-}\mathcal{X}$, and in particular definable subcategories of modules over rings. In the end, we give a conceptual proof of Auslander-Gruson-Jensen duality, which makes the duality between definable subcategories of left and right module more transparent.
	\end{abstract}
	
	\maketitle


	\section{Introduction}
	Let $\mathcal{X}$ be a skeletally small additive category.
	In this note, following Mitchell \cite{Mit72} we see $\mathcal{X}$ as a ring with several objects. A left $\mathcal{X}$-module is an additive covariant functor from $\mathcal{X}$ to $\Ab$, the category of all abelian groups.
	Similarly, a right $\mathcal{X}$-module is an additive contravariant functor from $\mathcal{X}$ to $\Ab$. We denote the category of all left (resp. right) $\mathcal{X}$-modules by $\Mod\text{-}\mathcal{X}$ (resp. $\Mod\text{-}\mathcal{X}^{op}$).
	Also the subcategory of finitely presented left (resp. right) $\mathcal{X}$-modules is denoted by $\modd\text{-}\mathcal{X}$ (resp. $\modd\text{-}\mathcal{X}^{op}$).
	
	Definable subcategories of $\Mod\text{-}\mathcal{X}$ play a central rule in model theory of modules \cite{Zie84,Pre09}. Definable subcategories can be defined using finitely presented functors on $\modd\text{-}\mathcal{X}$. A finitely presented functors on $\modd\text{-}\mathcal{X}$ is nothing but an object of the abelian category $\modd\text{-}(\modd\text{-}\mathcal{X})$, the {\it free abelian category} over $\mathcal{X}$.
	The free abelian category over $\mathcal{X}$ is, by definition, an abelian category $\mathcal{A}$ together with an additive functor $\mathcal{X}\rightarrow \mathcal{A}$, such that any other additive functor $\mathcal{X}\rightarrow \mathcal{B}$ to an abelian category $\mathcal{B}$ factors through $\mathcal{X}\rightarrow \mathcal{A}$ by a unique, up to natural isomorphism, exact functor $\mathcal{A}\rightarrow\mathcal{B}$ \cite{Fre65}. The existence of the free abelian category was proved by Freyd \cite{Fre65}. Constructing the free abelian category as the composition of two (contravariant) Yoneda embeddings
	\begin{center}
		$\mathcal{X}\rightarrow \modd\text{-}\mathcal{X}\rightarrow \modd\text{-}(\modd\text{-}\mathcal{X})$
	\end{center}
    is due to Gruson \cite{G}, Beligianis \cite[Theorem 6.1]{Bel00} and Krause \cite[Universal Property 2.10]{Kra98}.
	Before them, Adelman had given a simple construction of the free abelian category \cite{Ade73}.
	Objects of Adelman's free abelian category are diagrams (not necessarily complexes) of the form
	\begin{center}
		$X\rightarrow Y\rightarrow Z$
	\end{center}
	in $\mathcal{X}$, and morphisms are commutative diagrams modulo some homotopy relation (see the next section for details). We abuse the notation and denote by $\Romanbar{3}(\mathcal{X})$ Adelman's free abelian category over $\mathcal{X}$.
	
	Recall that a subcategory $\mathscr{D}$ of $\Mod\text{-}\mathcal{X}$ is called {\it definable} if it is closed under direct products, direct limits and pure subobjects. Another way of characterizing definable subcategories of $\Mod\text{-}\mathcal{X}$ is as follow: A subcategory $\mathscr{D}\subseteq\Mod\text{-}\mathcal{X}$ is definable if and only if there is a family $(\mathbb{F}_{\lambda})_{\lambda\in \Lambda}$ of functors in $\modd\text{-}(\modd\text{-}\mathcal{X})$, such that
	\begin{center}
		$\mathscr{D}=\{M\in \Mod\text{-}\mathcal{X}\mid \overrightarrow{\mathbb{F}_{\lambda}}(M)=0, \forall \lambda\in \Lambda\}$.
	\end{center}
	Where $\overrightarrow{\mathbb{F}_{\lambda}}$ is the unique extension of $\mathbb{F}_{\lambda}$ to $\Mod\text{-}\mathcal{X}$ which commutes with direct limits \cite[Corollary 10.2.32]{Pre09}. 
	
	By the universal property of free abelian category, there exists a unique canonical exact equivalence between $\Romanbar{3}(\mathcal{X})$ and $\modd$-$(\modd$-$\mathcal{X})$.
	The following question was the motivation for this work.
	\begin{question}
		Let $\mathcal{X}$ be a skeletally small additive category. Is it possible to characterize definable subcategories of $\Mod\text{-}\mathcal{X}$ using the category $\Romanbar{3}(\mathcal{X})$?
	\end{question}
    The main theorem, Theorem 3.9 will answer this question. Indeed we will prove that $\mathscr{D}\subseteq\Mod\text{-}\mathcal{X}$ is definable if and only if there is a family $(X_{\lambda}\overset{f_{\lambda}}{\rightarrow} Y_{\lambda}\overset{g_{\lambda}}{\rightarrow} Z_{\lambda})_{\lambda\in \Lambda}$ of objects in $\Romanbar{3}(\mathcal{X})$ such that
    \begin{align*}
    	\mathscr{D}=\{M\in \Mod\text{-}\mathcal{X} \mid \Ker(M(g_{\lambda}))\subseteq \Imm(M(f_{\lambda}))\; \forall \lambda\in \Lambda\}.
    \end{align*}
	
	Now let $A$ be a ring and $\proj(A)$ be the category of finitely generated projective left $A$-modules. Since $\Mod\text{-}A\cong\Mod\text{-}(\proj(A))$ and  $\modd\text{-}A\cong\modd\text{-}(\proj(A))$, we can apply our results to modules over $A$.
	In particular we prove that a subcategory $\mathscr{D}$ of $\Mod$-$A$ is definable if and only if there is a family of pairs of matrices $(U_{\lambda},V_{\lambda})$ (of appropriate size), with entries from $A$, such that 
	\begin{align*}
		\mathscr{D}=\{M\in \Mod\text{-}A \mid \forall \lambda\in \Lambda \;\text{and}\; \forall x\in M^n, U_{\lambda}x=0\implies x=V_{\lambda}y\; \text{for some}\; y\in M^m\}.
	\end{align*} 
	Auslander-Gruson-Jenson duality, stated in a general form as in \cite[Section 10.3]{Pre09}, states that there is a natural exact duality $d:\modd\text{-}(\modd\text{-}\mathcal{X}) \rightarrow \modd\text{-}(\modd\text{-}\mathcal{X}^{op})$. In the last section, we give a conceptual proof of this duality. Also using this, we show that there is a natural bijection between definable subcategories of $\Mod\text{-}\mathcal{X}$ and definable subcategories of $\Mod\text{-}\mathcal{X}^{op}$. This result for module category of rings is due to Herzog \cite{Her93} (see \cite{Baz08}).

	
	\section{preliminaries}
	In this section we recall the definition of modules over an additive category following Mitchell \cite{Mit72}. Then we give two presentations of the free abelian category and some of their basic properties that we will need in the sequel.
	
	Let $\mathcal{X}$ be a skeletally small additive category. By a {\it left $\mathcal{X}$-module} we mean a covariant additive functor
	\begin{center}
		$M:\mathcal{X}\longrightarrow \Ab$,
	\end{center}
	where $\Ab$ is the category of all abelian groups. The collection of all left $\mathcal{X}$-modules together with natural transformations between them forms an abelian category that is denoted by $\Mod\text{-}\mathcal{X}$. For every $X\in\mathcal{X}$ the representable functor $\Hom_{\mathcal{X}}(X,-)$ is a projective left $\mathcal{X}$-module and by the Yoneda lemma $\Hom_{\mathcal{X}}(X,Y)=\Hom_{\Mod\text{-}\mathcal{X}}(\Hom_{\mathcal{X}}(Y,-),\Hom_{\mathcal{X}}(X,-))$. A left $\mathcal{X}$-module $M$ is called {\it finitely presented} if there exists an exact sequence
	\begin{center}
		$\Hom_{\mathcal{X}}(Y,-)\rightarrow \Hom_{\mathcal{X}}(X,-)\rightarrow M\rightarrow 0$.
	\end{center}
    in $\Mod\text{-}\mathcal{X}$.
    The subcategory of all finitely presented left $\mathcal{X}$-modules is denoted by $\modd\text{-}\mathcal{X}$. Similarly a right $\mathcal{X}$-module is a contravariant functor from $\mathcal{X}$ to $\Ab$, or equivalently a covariant functor $N:\mathcal{X}^{op}\longrightarrow \Ab$. The category of all (resp. finitely presented) right $\mathcal{X}$-modules is denoted by $\Mod\text{-}\mathcal{X}^{op}$ (resp. $\modd\text{-}\mathcal{X}^{op}$).
    
    \begin{remark}
    	In \cite{Mit72} (and almost all other sources) modules over pre-additive categories are in consideration.
    	If $\mathcal{X}$ is a pre-additive category, then $\Mod\text{-}\mathcal{X}\cong \Mod\text{-}(\proj(\mathcal{X}))$, where $\proj(\mathcal{X})$ is the category of all finitely generated projective objects in $\Mod\text{-}\mathcal{X}$. So we don't lose anything if we only consider modules over additive categories.
    \end{remark}
   In the following we recall the definition of the free abelian category over an additive category, and also we give two different presentations of the free abelian category, which are our main tools for characterizing definable subcategories.
	\begin{definition}\cite{Fre65}
		Let $\mathcal{X}$ be an additive category. The free abelian category over $\mathcal{X}$ is by definition an abelian category $\mathcal{A}$ together with an additive functor $\mathcal{X}\rightarrow \mathcal{A}$ such that, for any other abelian category $\mathcal{B}$ and any additive functor $\mathcal{X}\rightarrow \mathcal{B}$, there exists a unique, up to natural isomorphism, exact functor $\mathcal{A}\rightarrow \mathcal{B}$ making the following diagram commutative.
	    
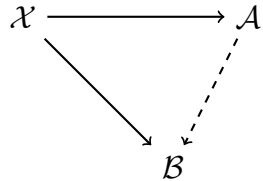
\begin{center}
	    	\begin{tikzpicture}
	    		\node (X1) at (0,0) {$\mathcal{X}$};
	    		\node (X2) at (3,0) {$\mathcal{A}$};
	    		\node (X3) at (2,-2) {$\mathcal{B}$};
	    		\draw [->,thick] (X1) -- (X2) node [midway,above] {};
	    		\draw [->,thick] (X1) -- (X3) node [midway,above] {};
	    		\draw [->,thick,dashed] (X2) -- (X3) node [midway,left] {};
	    	\end{tikzpicture}
	    \end{center}	
	\end{definition}
    
    The free abelian category, if it exists, is unique up to a unique natural isomorphism due to its universal property. In \cite{Fre65}, Freyd proved that the free abelian category always exists. In what follows we will give two different presentations of the free abelian category.
    
    Let $\mathcal{X}$ be a skeletally small additive category and $\rm Ch^{1,2,3}(\mathcal{X})$ be the category of three term chains (not necessarily complex) over $\mathcal{X}$. We denote a test object of $\rm Ch^{1,2,3}(\mathcal{X})$ by $X_1\overset{f_1}{\rightarrow} X_2\overset{f_2}{\rightarrow} X_3$ and a morphism by a commutative diagram like
    \begin{equation}\tag{Diagram 2.1}
    	\begin{tikzpicture}\label{dia2.1}
    		\node (X1) at (-2,0) {$X_1$};
    		\node (X2) at (0,0) {$X_2$};
    		\node (X3) at (2,0) {$X_3$};
    		\node (X4) at (-2,-2) {$Y_1$};
    		\node (X5) at (0,-2) {$Y_2$};
    		\node (X6) at (2,-2) {$Y_3$};
    		\draw [->,thick] (X1) -- (X2) node [midway,above] {$f_1$};
    		\draw [->,thick] (X2) -- (X3) node [midway,above] {$f_2$};
    		\draw [->,thick] (X4) -- (X5) node [midway,above] {$g_1$};
    		\draw [->,thick] (X5) -- (X6) node [midway,above] {$g_2$};
    		\draw [->,thick] (X1) -- (X4) node [midway,left] {$\alpha_1$};
    		\draw [->,thick] (X2) -- (X5) node [midway,left] {$\alpha_2$};
    		\draw [->,thick] (X3) -- (X6) node [midway,left] {$\alpha_3$};
    	\end{tikzpicture}
    \end{equation}
    or simply by the triple $(\alpha_1,\alpha_2,\alpha_3)$. We say that $(\alpha_1,\alpha_2,\alpha_3)$ is null-homotopic if there are morphisms $s:X_2\rightarrow Y_1$ and $t:X_3\rightarrow Y_2$ such that $\alpha_2= g_1s+tf_2$. It is not hard to see that null-homotopic morphisms form an ideal of the additive category $\rm Ch^{1,2,3}(\mathcal{X})$. We denote by $\Romanbar{3}(\mathcal{X})$ the associated factor category.
    
    \begin{theorem}\cite{Ade73}\label{2.3}
    	Let $\mathcal{X}$ be a skeletally small additive category. Then $\Romanbar{3}(\mathcal{X})$ is an abelian category. Moreover the canonical full embedding $\mathfrak{j}:\mathcal{X}\rightarrow \Romanbar{3}(\mathcal{X})$, given by $X\mapsto (0\rightarrow X\rightarrow 0)$, provides the free abelian category over $\mathcal{X}$.
    \end{theorem}
    By the canonical embedding $\mathfrak{j}:\mathcal{X}\rightarrow \Romanbar{3}(\mathcal{X})$, mentioned in the above theorem, sometimes we identify $X\in \mathcal{X}$ with its image $\mathfrak{j}(X)=(0\rightarrow X\rightarrow 0)$.
    
    For the future purposes, we want to tell something about the abelian structure of $\Romanbar{3}(\mathcal{X})$. The proofs can be found in \cite{Ade73}. Let $\alpha=(\alpha_1,\alpha_2,\alpha_3)$ be a morphism in $\Romanbar{3}(\mathcal{X})$ as in the \ref{dia2.1}. Then the kernel and cokernel of $\alpha$ are as follows.
    \begin{center}
    	\begin{tikzpicture}
    		\node (K0) at (-3,3) {$\Ker(\alpha)$};
    		\node (K1) at (0,3) {$X_1\oplus Y_1$};
    		\node (K2) at (4,3) {$X_2\oplus Y_1$};
    		\node (K3) at (8,3) {$X_3\oplus Y_2$};
    		\node (X0) at (-3,0) {$X$};
    		\node (X1) at (0,0) {$X_1$};
    		\node (X2) at (4,0) {$X_2$};
    		\node (X3) at (8,0) {$X_3$};
    		\node (Y0) at (-3,-3) {$Y$};
    		\node (Y1) at (0,-3) {$Y_1$};
    		\node (Y2) at (4,-3) {$Y_2$};
    		\node (Y3) at (8,-3) {$Y_3$};
    		\node (C0) at (-3,-6) {$\Coker(\alpha)$};
    		\node (C1) at (0,-6) {$Y_1\oplus X_2$};
    		\node (C2) at (4,-6) {$Y_2\oplus X_3$};
    		\node (C3) at (8,-6) {$Y_3\oplus X_3$};
    		\draw [->,thick] (K1) -- (K2) node [midway,above] {$\begin{bmatrix}f_1&0\\ \alpha_1&-1\end{bmatrix}$};
    		\draw [->,thick] (K2) -- (K3) node [midway,above] {$\begin{bmatrix}f_2&0\\\alpha_2&-g_1\end{bmatrix}$};
    		\draw [->,thick] (X1) -- (X2) node [midway,above] {$f_1$};
    		\draw [->,thick] (X2) -- (X3) node [midway,above] {$f_2$};
    		\draw [->,thick] (Y1) -- (Y2) node [midway,above] {$g_1$};
    		\draw [->,thick] (Y2) -- (Y3) node [midway,above] {$g_2$};
    		\draw [->,thick] (C1) -- (C2) node [midway,above] {$\begin{bmatrix}g_1&\alpha_2\\0&-f_2\end{bmatrix}$};
    		\draw [->,thick] (C2) -- (C3) node [midway,above] {$\begin{bmatrix}g_2&\alpha_3\\0&-1\end{bmatrix}$};
    		\draw [->,thick] (K0) -- (X0) node [midway,left] {};
    		\draw [->,thick] (X0) -- (Y0) node [midway,left] {$\alpha$};
    		\draw [->,thick] (Y0) -- (C0) node [midway,above] {};
    		\draw [->,thick] (K1) -- (X1) node [midway,left] {$\begin{bmatrix}1&0\end{bmatrix}$};
        	\draw [->,thick] (X1) -- (Y1) node [midway,left] {$\alpha_1$};
    		\draw [->,thick] (Y1) -- (C1) node [midway,left] {$\begin{bmatrix}1\\0\end{bmatrix}$};
    		\draw [->,thick] (K2) -- (X2) node [midway,left] {$\begin{bmatrix}1&0\end{bmatrix}$};
    		\draw [->,thick] (X2) -- (Y2) node [midway,left] {$\alpha_2$};
    		\draw [->,thick] (Y2) -- (C2) node [midway,left] {$\begin{bmatrix}1\\0\end{bmatrix}$};
    		\draw [->,thick] (K3) -- (X3) node [midway,left] {$\begin{bmatrix}1&0\end{bmatrix}$};
    		\draw [->,thick] (X3) -- (Y3) node [midway,left] {$\alpha_3$};
    		\draw [->,thick] (Y3) -- (C3) node [midway,left] {$\begin{bmatrix}1\\0\end{bmatrix}$};
    	\end{tikzpicture}
    \end{center}
    In particular, $0\rightarrow Y\rightarrow Z$ is the kernel of $\mathfrak{j}(Y)\rightarrow \mathfrak{j}(Z)$ and $X\rightarrow Y\rightarrow 0$ is the cokernel of $\mathfrak{j}(X)\rightarrow \mathfrak{j}(Y)$.
    
    \begin{proposition}\cite[Proposition 1.5]{Ade73}\label{2.4}
    	Let $X\overset{f}{\rightarrow}Y\overset{g}{\rightarrow}Z$ be an object in $\Romanbar{3}(\mathcal{X})$. Then we have a commutative diagram\begin{center}
    		\begin{tikzpicture}
    			\node (X1) at (-6,0) {$\mathfrak{j}(X)$};
    			\node (X2) at (0,0) {$\mathfrak{j}(Y)$};
    			\node (X3) at (6,0) {$\mathfrak{j}(Z)$};
    			\node (X4) at (-3,-3) {$(0\rightarrow Y\rightarrow Z)$};
    			\node (X5) at (3,-3) {$(X\rightarrow Y\rightarrow 0)$};
    			\draw [->,thick] (X1) -- (X2) node [midway,above] {$\mathfrak{j}(f)$};
    			\draw [->,thick] (X2) -- (X3) node [midway,above] {$\mathfrak{j}(g)$};
    			\draw [->,thick] (X4) -- (X2) node [midway,above] {$k=(0,1,0)$};
    			\draw [->,thick] (X2) -- (X5) node [midway,above] {$c=(0,1,0)$};
    			\draw [->,thick] (X4) -- (X5) node [midway,above] {$m=(0,1,0)$};
    		\end{tikzpicture}
    	\end{center}
     in $\Romanbar{3}(\mathcal{X})$ such that $k$ is the kernel of $\mathfrak{j}(g)$, $c$ is the cokernel of $\mathfrak{j}(f)$ and $X\overset{f}{\rightarrow}Y\overset{g}{\rightarrow}Z$ is isomorphic to $\Imm(m)$.
    \end{proposition}

    For the other presentation of the free abelian category, note that $\modd\text{-}\mathcal{X}$ is an additive category with cokernels. So by a result of Freyd \cite{Fre65}, $\modd\text{-}(\modd\text{-}\mathcal{X})$ is an abelian category. Moreover:
    
    \begin{theorem}\cite{G,Bel00,Kra98}\label{2.5}
    	The composition of two Yoneda embeddings
    	\begin{center}
    		$\mathfrak{i}:\mathcal{X}\rightarrow \modd\text{-}\mathcal{X} \rightarrow \modd\text{-}(\modd\text{-}\mathcal{X})$
    	\end{center}
    	provides the free abelian category over $\mathcal{X}$.
    \end{theorem}

    \begin{remark}\label{2.6}
    	Let $\mathbb{F}$ be a functor in  $\modd\text{-}(\modd\text{-}\mathcal{X})$. So there is an exact sequence
    	\begin{center}
    		$\Hom(F_2,-)\rightarrow\Hom(F_1,-)\rightarrow \mathbb{F}\rightarrow 0$.
    	\end{center}
        In $\Mod\text{-}(\modd\text{-}\mathcal{X})$ with $F_1,F_2\in \modd\text{-}\mathcal{X}$. Then we are given the following exact sequences in $\Mod\text{-}\mathcal{X}$
        \begin{align*}
        	\Hom(Y_1,-)\rightarrow\Hom(X_1,-)\rightarrow F_1\rightarrow 0,\\
        	\Hom(Y_2,-)\rightarrow\Hom(X_2,-)\rightarrow F_2\rightarrow 0.
        \end{align*}
        Altogether, we have the following commutative diagram with exact row and exact columns in $\Mod\text{-}(\modd\text{-}\mathcal{X})$.
    \end{remark}
\begin{equation}\tag{Diagram 2.2}
	\begin{tikzpicture}\label{dia2.2}
		\node (K2) at (-3,6) {$0$};
		\node (K1) at (0,6) {$0$};
		\node (F2) at (-3,3) {$\bigl(F_2,-\bigl)$};
		\node (F1) at (0,3) {$\bigl(F_1,-\bigl)$};
		\node (FF) at (3,3) {$\mathbb{F}$};
		\node (F0) at (6,3) {$0$};
		\node (X2) at (-3,0) {$\bigl((X_2,-),-\bigl)$};
		\node (X1) at (0,0) {$\bigl((X_1,-),-\bigl)$};
		\node (Y2) at (-3,-3) {$\bigl((Y_2,-),-\bigl)$};
		\node (Y1) at (0,-3) {$\bigl((Y_1,-),-\bigl)$};
		\draw [->,thick] (K2) -- (F2) node [midway,above] {};
		\draw [->,thick] (K1) -- (F1) node [midway,above] {};
		\draw [->,thick] (F2) -- (X2) node [midway,above] {};
		\draw [->,thick] (X2) -- (Y2) node [midway,above] {};
		\draw [->,thick] (F1) -- (X1) node [midway,above] {};
		\draw [->,thick] (X1) -- (Y1) node [midway,above] {};
		\draw [->,thick] (F2) -- (F1) node [midway,above] {};
		\draw [->,thick] (F1) -- (FF) node [midway,above] {};
		\draw [->,thick] (FF) -- (F0) node [midway,left] {};
		\draw [->,thick] (X2) -- (X1) node [midway,left] {};
		\draw [->,thick] (Y2) -- (Y1) node [midway,above] {};
	\end{tikzpicture}
\end{equation}
     By Yoneda lemma, the lower square is induced by a commutative square
     \begin{equation}\tag{Diagram 2.3}
     	\begin{tikzpicture}\label{dia2.3}
     		\node (X1) at (-3,1) {$X_2$};
     		\node (X2) at (0,1) {$X_1$};
     		\node (X3) at (-3,-1) {$Y_2$};
     		\node (X4) at (0,-1) {$Y_1$};
     		\draw [->,thick] (X1) -- (X2) node [midway,above] {$f$};
     		\draw [->,thick] (X1) -- (X3) node [midway,left] {$a$};
     		\draw [->,thick] (X2) -- (X4) node [midway,left] {$b$};
     		\draw [->,thick] (X3) -- (X4) node [midway,above] {$g$};
     	\end{tikzpicture}
     \end{equation}
     in $\mathcal{X}$.
 
    \begin{theorem}\label{2.7}
    	Let $\mathcal{X}$ be an additive category. Then, up to natural isomorphism, there is a unique exact functor $\mathfrak{k}:\modd$-$(\modd$-$\mathcal{X})\rightarrow \Romanbar{3}(\mathcal{X})$ which makes the following diagram commutative.
    	\begin{center}
    		\begin{tikzpicture}
    			\node (X1) at (0,0) {$\mathcal{X}$};
    			\node (X2) at (3,0) {$\modd$-$(\modd$-$\mathcal{X})$};
    			\node (X3) at (2,-2) {$\Romanbar{3}(\mathcal{X})$};
    			\draw [->,thick] (X1) -- (X2) node [midway,above] {$\mathfrak{i}$};
    			\draw [->,thick] (X1) -- (X3) node [midway,left] {$\mathfrak{j}$};
    			\draw [->,thick,dashed] (X2) -- (X3) node [midway,right] {$\mathfrak{k}$};
    		\end{tikzpicture}
    	\end{center}
    	For an arbitrary object $\mathbb{F}\in \modd$-$(\modd$-$\mathcal{X})$, if we construct \ref{dia2.2} and \ref{dia2.3} for $\mathbb{F}$, then 
    	\begin{center}
    		$\mathfrak{k}(\mathbb{F})\cong X_2\overset{\begin{bmatrix}-a\\f\end{bmatrix}}{\longrightarrow}X_1\oplus Y_2\overset{\begin{bmatrix}b&g\\0&-1\end{bmatrix}}{\longrightarrow}Y_1\oplus Y_2$
    	\end{center}
    	\begin{proof}
    		By the definition of free abelian category, the exact functor $\mathfrak{k}$ exists. First note that $\mathfrak{k}(\bigl((X,-),-\bigl))=0\rightarrow X\rightarrow 0$.
    		Now let $\mathbb{F}$ be an arbitrary object in $\modd$-$(\modd$-$\mathcal{X})$ and suppose that \ref{dia2.2} and \ref{dia2.3} are given. By exactness of $\mathfrak{k}$ we have that 
    		\begin{align*}
    			\mathfrak{k}(F_1)=0\rightarrow X_1\overset{b}{\rightarrow} Y_1,\\
    			\mathfrak{k}(F_2)=0\rightarrow X_2\overset{a}{\rightarrow} Y_2.
    		\end{align*}
    		Again by exactness of $\mathfrak{k}$, $\mathfrak{k}(\mathbb{F})$ is isomorphic to the cokernel of the following morphism in $\Romanbar{3}(\mathcal{X})$.
    		\begin{center}
    			\begin{tikzpicture}
    				\node (X1) at (-2,0) {$0$};
    				\node (X2) at (0,0) {$X_2$};
    				\node (X3) at (2,0) {$Y_2$};
    				\node (X4) at (-2,-2) {$0$};
    				\node (X5) at (0,-2) {$X_1$};
    				\node (X6) at (2,-2) {$Y_1$};
    				\draw [->,thick] (X1) -- (X2) node [midway,above] {};
    				\draw [->,thick] (X2) -- (X3) node [midway,above] {$a$};
    				\draw [->,thick] (X4) -- (X5) node [midway,above] {};
    				\draw [->,thick] (X5) -- (X6) node [midway,above] {$b$};
    				\draw [->,thick] (X1) -- (X4) node [midway,left] {};
    				\draw [->,thick] (X2) -- (X5) node [midway,left] {$f$};
    				\draw [->,thick] (X3) -- (X6) node [midway,left] {$g$};
    			\end{tikzpicture}
    		\end{center}
    		Then the result follows from the construction of cokernel in $\Romanbar{3}(\mathcal{X})$.
    	\end{proof}
    \end{theorem}
 
    For the quasi inverse of the equivalence $\mathfrak{k}$ we need a version of snake lemma. 
    \begin{lemma}\label{2.8}
    	Let $X\overset{f}{\rightarrow}Y\overset{g}{\rightarrow}Z$ be a adiagram in any abelian category. Then we have the following exact sequence.
    	\begin{align*}
    		0\rightarrow &\Ker(f)\rightarrow\ker(gf)\rightarrow Ker(g)\rightarrow\\ &\Coker(f)\rightarrow\Coker(gf)\rightarrow Coker(g)\rightarrow 0.
    	\end{align*}
    \end{lemma}
  
	Now we can construct the quasi inverse of $\mathfrak{k}$.
	\begin{theorem}\label{2.9}
		Let $\mathcal{X}$ be an additive category. Then there is a quasi inverse $\mathfrak{k}^{-1}:\Romanbar{3}(\mathcal{X})\rightarrow \modd$-$(\modd$-$\mathcal{X})$ for $\mathfrak{k}$, which makes the diagram
		\begin{center}
			\begin{tikzpicture}
				\node (X1) at (0,0) {$\mathcal{X}$};
				\node (X2) at (3,0) {$\modd$-$(\modd$-$\mathcal{X})$};
				\node (X3) at (2,-2) {$\Romanbar{3}(\mathcal{X})$};
				\draw [->,thick] (X1) -- (X2) node [midway,above] {$\mathfrak{i}$};
				\draw [->,thick] (X1) -- (X3) node [midway,left] {$\mathfrak{j}$};
				\draw [->,thick,dashed] (X3) -- (X2) node [midway,right] {$\mathfrak{k}^{-1}$};
			\end{tikzpicture}
		\end{center}
	    commutative, and for every object $X\overset{f}{\rightarrow}Y\overset{g}{\rightarrow}Z\in \Romanbar{3}(\mathcal{X})$ we have
		\begin{center}
			$\mathfrak{k}^{-1}(X\overset{f}{\rightarrow}Y\overset{g}{\rightarrow}Z)\cong \mathbb{F}$, 
		\end{center}
		where $\mathbb{F}$ placed in the following exact diagram.
		\begin{center}
			\begin{tikzpicture}
				\node (K2) at (-3,6) {$0$};
				\node (K1) at (0,6) {$0$};
				\node (F2) at (-3,3) {$\bigl(G,-\bigl)$};
				\node (F1) at (0,3) {$\bigl(F,-\bigl)$};
				\node (FF) at (3,3) {$\mathbb{F}$};
				\node (F0) at (6,3) {$0$};
				\node (X2) at (-3,0) {$\bigl((X,-),-\bigl)$};
				\node (X1) at (0,0) {$\bigl((Y,-),-\bigl)$};
				\node (Y2) at (-3,-3) {$\bigl((Z,-),-\bigl)$};
				\node (Y1) at (0,-3) {$\bigl((Z,-),-\bigl)$};
				\draw [->,thick] (K2) -- (F2) node [midway,above] {};
				\draw [->,thick] (K1) -- (F1) node [midway,above] {};
				\draw [->,thick] (F2) -- (X2) node [midway,above] {};
				\draw [->,thick] (X2) -- (Y2) node [midway,left] {$\tilde{gf}$};
				\draw [->,thick] (F1) -- (X1) node [midway,above] {};
				\draw [->,thick] (X1) -- (Y1) node [midway,right] {$\tilde{g}$};
				\draw [->,thick] (F2) -- (F1) node [midway,above] {};
				\draw [->,thick] (F1) -- (FF) node [midway,above] {};
				\draw [->,thick] (FF) -- (F0) node [midway,left] {};
				\draw [->,thick] (X2) -- (X1) node [midway,above] {$\tilde{f}$};
				\draw [double,-,thick] (Y2) -- (Y1) node [midway,above] {};
			\end{tikzpicture}
		\end{center}
		
		\begin{proof}
			By the definition of free abelian category, the exact functor $\mathfrak{k}^{-1}$ exists. First note that $\mathfrak{k}^{-1}(0\rightarrow X\rightarrow 0)=\bigl((X,-),-\bigl)$.
			Now let $\mathbb{X}=(X\overset{f}{\rightarrow}Y\overset{g}{\rightarrow}Z)$ be an arbitrary object in $\Romanbar{3}(\mathcal{X})$. By Proposition \ref{2.4} and Lemma \ref{2.8} we have the following natural isomorphisms
			\begin{align*}
				\mathbb{X}&\cong \Imm\big(\Ker(\mathfrak{j}(g))\rightarrow \Coker(\mathfrak{j}(f))\big)\\
				          &\cong \Coker\big(\Ker(\mathfrak{j}(gf))\rightarrow \Ker(\mathfrak{j}(g))\big).
			\end{align*}
		   Then by exactness of $\mathfrak{k}^{-1}$,
		   \begin{align*}
		   	&\mathfrak{k}^{-1}\big(\Ker(\mathfrak{j}(gf))\big)\cong \big(G,-\big),\\
		   	&\mathfrak{k}^{-1}\big(\Ker(\mathfrak{j}(g))\big)\cong \big(F,-\big)\\
		   	&\text{and so}\\
		   	&\mathfrak{k}^{-1}(\mathbb{X})\cong \mathbb{F}.
		   \end{align*}
		\end{proof}
	\end{theorem}
	\section{Characterization of definable subcategories}
	In this section, using the results of the previous section, we give a new characterization of definable subcategories. First let us recall some information about definable subcategories.
	
	Let $\mathcal{X}$ be a skeletally small additive category. We are interested in definable subcategories of $\Mod\text{-}\mathcal{X}$ and $\Mod\text{-}\mathcal{X}^{op}$.
	\begin{definition}
		An additive functor $M:\Mod\text{-}\mathcal{X}\rightarrow \Ab$ is called {\it coherent} if it satisfies one of the following two equivalent conditions.
		\begin{itemize}
			\item[(1)]
			$M$ commutes with direct limits, and it is finitely presented when restricted to $\modd\text{-}\mathcal{X}$.
			\item[(2)]
			There is a sequence of natural transformations
			\begin{center}
				$(G,-)\rightarrow (F,-)\rightarrow M\rightarrow 0$,
			\end{center}
		    with $F,G\in \modd\text{-}\mathcal{X}$, which is exact when evaluated at any object in $\Mod\text{-}\mathcal{X}$.
		\end{itemize}
	\end{definition}

    \begin{remark}\label{3.2}
    	Any functor $\mathbb{F}:\modd\text{-}\mathcal{X}\rightarrow \Ab$ can uniquely extend to a functor $\overrightarrow{\mathbb{F}}:\Mod\text{-}\mathcal{X} \rightarrow \Ab$ which commutes with direct limits. So we can say that coherent functors on $\Mod\text{-}\mathcal{X}$ are exactly the unique extension of finitely presented functors from $\modd\text{-}\mathcal{X}$ to $\Ab$.
    \end{remark}

    Now we can define definable subcategories. 
    
    \begin{definition}
    	Let $\mathcal{X}$ be a skeletally small additive category. A subcategory $\mathscr{D}\subseteq \Mod\text{-}\mathcal{X}$ is called {\it definable} if it satisfies one of the following equivalent conditions.
    		\begin{itemize}
    		\item[(1)] $\mathscr{D}$ is closed under direct products, direct limits and pure subobjects.
    		\item[(2)] There is a family $(\overrightarrow{\mathbb{F}}_{\lambda})_{\lambda}$ of coherent functors from $\Mod\text{-}\mathcal{X}$ to $\Ab$ such that 
    		\begin{align*}
    			\mathscr{D}=\bigcap_{\lambda \in \Lambda}\Ker(\overrightarrow{\mathbb{F}}_{\lambda}):=\{M\in \Mod\text{-}\mathcal{X} \mid \overrightarrow{\mathbb{F}}_{\lambda}(M)=0, \; \forall \lambda\in\Lambda\}.
    		\end{align*}
    	\end{itemize}
    \end{definition}
    For the equivalence of conditions in the above definition see for example
    \cite[Corollary 10.2.32.]{Pre09}. We want to give another algebraic characterization of definable subcaategories, which is a little simpler, and useful.
    
    In Section 2, we construct an equivalence of categories
    \begin{center}
    	\begin{tikzpicture}
    		\node (X1) at (0,0) {$\modd$-$(\modd$-$\mathcal{X})$};
    		\node (X2) at (4,0) {$\Romanbar{3}(\mathcal{X})$};
    		\node (k) at (2,1.3) {$\mathfrak{k}$};
    		\node (k') at (2,-0.7) {$\mathfrak{k}^{-1}$};
    		\draw [->,thick] (X1) to [out=45,in=135] (X2) node [midway,above] {};
    		\draw [->,thick] (X2) to [out=225,in=315] (X1) node [midway,above] {};
    	\end{tikzpicture}
    \end{center}
    By Remark \ref{3.2}, a definable subcategory of $\Mod\text{-}\mathcal{X}$ is nothing but the kernel of a family of functors of the form $\overrightarrow{\mathbb{F}}_{\lambda}$, with $\mathbb{F}\in \modd$-$(\modd$-$\mathcal{X})$. From the canonical equivalence $\mathfrak{k}$ it is natural to ask, what is the counterpart property of "$M\in \Ker(\overrightarrow{\mathbb{F}})$" in $\Romanbar{3}(\mathcal{X})$? In the sequel, we answer this question and using it we provide another characterization of definable subcategories.

	\begin{lemma}\label{3.4}
		Let $\mathcal{A}$ be an abelian category and
		\begin{center}
			\begin{tikzpicture}
				\node (X1) at (0,3) {$Q_1$};
				\node (X2) at (2,3) {$Q_2$};
				\node (X3) at (0,1.5) {$P_1$};
				\node (X4) at (2,1.5) {$P_2$};
				\node (X5) at (0,0) {$F_1$};
				\node (X6) at (2,0) {$F_2$};
				\node (X7) at (0,-1.5) {$0$};
				\node (X8) at (2,-1.5) {$0$};
				\draw [->,thick] (X1) -- (X2) node [midway,above] {$b$};
				\draw [->,thick] (X3) -- (X4) node [midway,above] {$c$};
				\draw [->,thick] (X5) -- (X6) node [midway,above] {};
				\draw [->,thick] (X1) -- (X3) node [midway,left] {$a$};
				\draw [->,thick] (X3) -- (X5) node [midway,left] {};
				\draw [->,thick] (X5) -- (X7) node [midway,left] {};
				\draw [->,thick] (X2) -- (X4) node [midway,left] {$d$};
				\draw [->,thick] (X4) -- (X6) node [midway,left] {};
				\draw [->,thick] (X6) -- (X8) node [midway,left] {};
			\end{tikzpicture}
		\end{center}
		be a commutative diagram with exact columns in $\mathcal{A}$. Then the following conditions are equivalent for any object $M\in \mathcal{A}$.
		\begin{itemize}
			\item[(1)]
			The induced map $\Hom(F_2,M)\rightarrow \Hom(F_1,M)$ is surjective.
			\item[(2)]
			In the diagram 
			\begin{center}
				\begin{tikzpicture}
					\node (X1) at (-4,0) {$Q_1\oplus Q_2$};
					\node (X2) at (0,0) {$P_1\oplus Q_2$};
					\node (X3) at (4,0) {$P_2$};
					\node (X4) at (0,-3) {$M$};
					\draw [->,thick] (X1) -- (X2) node [midway,above] {$\begin{bmatrix}a&0\\b&-1\end{bmatrix}$};
					\draw [->,thick] (X2) -- (X3) node [midway,above] {$\begin{bmatrix}c&-d\end{bmatrix}$};
					\draw [->,thick] (X2) -- (X4) node [midway,left] {$f$};
					\draw [->,thick,dashed] (X3) -- (X4) node [midway,right] {$g$};
					\draw [->,thick] (X1) -- (X4) node [midway,left] {$0$};
				\end{tikzpicture}
			\end{center}
		  for any morphism $f$ that makes the left-hand triangle commutative, there exists a morphism $g$ that makes the right-hand triangle commutative.
		\end{itemize}
	     \begin{proof}
	     	The proof is easy and it is left to the reader. Note that if $f=[f_1\quad f_2]$ make the left-hand triangle commutative, then $f_2=0$. And a morphism $F_i\rightarrow M$ is nothing but a morphism $P_i\rightarrow M$ with $Q_i\rightarrow P_i\rightarrow M=0$.
	     \end{proof}
	\end{lemma}

    \begin{lemma}\label{3.5}
    	Let $\mathbb{F}\in \modd\text{-}(\modd\text{-}\mathcal{X})$, and $X\overset{f}{\rightarrow} Y\overset{g}{\rightarrow}Z=\mathfrak{k}(\mathbb{F})$ be the corresponding object in $\Romanbar{3}(\mathcal{X})$ under the canonical equivalence $\mathfrak{k}$. Then the following are equivalent for an object $M\in \Mod\text{-}\mathcal{X}$.
    	\begin{itemize}
    		\item[(1)] $M\in \Ker(\overrightarrow{\mathbb{F}})$.
    		\item[(2)] $\Ker(M(g))\subseteq \Imm(M(f))$.
    	\end{itemize}
    \begin{proof}
    	By Theorem \ref{2.7}, if we consider the \ref{dia2.2} and \ref{dia2.3} for $\mathbb{F}$, then 
    		\begin{center}
    		$\mathfrak{k}(\mathbb{F})\cong X_2\overset{\begin{bmatrix}-a\\f\end{bmatrix}}{\longrightarrow}X_1\oplus Y_2\overset{\begin{bmatrix}b&g\\0&-1\end{bmatrix}}{\longrightarrow}Y_1\oplus Y_2$
    	\end{center}
        From the \ref{dia2.2} we have that $M\in \Ker(\overrightarrow{\mathbb{F}})$ if and only if the induced map $\Hom(F_2,M)\rightarrow \Hom(F_1,M)$ is surjective. On the other hand, \ref{dia2.2} gives the following commutative diagram with exact columns.
        \begin{center}
        	\begin{tikzpicture}
        		\node (X1) at (0,3) {$(Y_1,-)$};
        		\node (X2) at (2,3) {$(Y_2,-)$};
        		\node (X3) at (0,1.5) {$(X_1,-)$};
        		\node (X4) at (2,1.5) {$(X_2,-)$};
        		\node (X5) at (0,0) {$F_1$};
        		\node (X6) at (2,0) {$F_2$};
        		\node (X7) at (0,-1.5) {$0$};
        		\node (X8) at (2,-1.5) {$0$};
        		\draw [->,thick] (X1) -- (X2) node [midway,above] {};
        		\draw [->,thick] (X3) -- (X4) node [midway,above] {};
        		\draw [->,thick] (X5) -- (X6) node [midway,above] {};
        		\draw [->,thick] (X1) -- (X3) node [midway,left] {};
        		\draw [->,thick] (X3) -- (X5) node [midway,left] {};
        		\draw [->,thick] (X5) -- (X7) node [midway,left] {};
        		\draw [->,thick] (X2) -- (X4) node [midway,left] {};
        		\draw [->,thick] (X4) -- (X6) node [midway,left] {};
        		\draw [->,thick] (X6) -- (X8) node [midway,left] {};
        	\end{tikzpicture}
        \end{center}
    Now the result follow from Lemma \ref{3.4} and Yoneda lemma.
    \end{proof}
    \end{lemma}

    \begin{lemma}\label{3.6}
    	Let $\mathbb{X}=X\overset{f}{\rightarrow} Y\overset{g}{\rightarrow}Z\in \Romanbar{3}(\mathcal{X})$ and
    	$\mathbb{F}=\mathfrak{k}^{-1}(\mathbb{X})$ be the corresponding object in $\modd\text{-}(\modd\text{-}\mathcal{X})$ under the canonical equivalence $\mathfrak{k}^{-1}$.
    	Then the following are equivalent for an object $M\in \Mod\text{-}\mathcal{X}$.
    	\begin{itemize}
    		\item[(1)] $M\in \Ker(\overrightarrow{\mathbb{F}})$.
    		\item[(2)] $\Ker(M(g))\subseteq \Imm(M(f))$.
    	\end{itemize}
        \begin{proof}
        	The proof uses Theorem \ref{2.9} and it is very similar to the proof of Lemma \ref{3.5}.
        \end{proof}
    \end{lemma}

    Using the above result we can prove the main theorem, which gives a new characterization of definable subcategories. First let fix a notation to express results more clearly.
    
    \begin{notation}\label{3.7}
    	Let $\mathbb{X}=X\overset{f}{\rightarrow} Y\overset{g}{\rightarrow}Z$ be a diagram in additive category $\mathcal{X}$. Define
    	\begin{center}
    		$\Omega_{\mathbb{X}}:=\{M\in \Mod\text{-}\mathcal{X} \mid \Ker(M(g))\subseteq \Imm(M(f))\}$.
    	\end{center}
    \end{notation}

    \begin{remark}\label{3.8}
    	\begin{itemize}
    		\item[(1)]
    		By Yoneda lemma, $\Omega_{\mathbb{X}}$ is equal to the class of all functors $M\in \Mod\text{-}\mathcal{X}$ such that in the diagram 
    		\begin{center}
    			\begin{tikzpicture}
    				\node (X1) at (-4,0) {$(Z,-)$};
    				\node (X2) at (0,0) {$(Y,-)$};
    				\node (X3) at (4,0) {$(X,-)$};
    				\node (X4) at (0,-3) {$M$};
    				\draw [->,thick] (X1) -- (X2) node [midway,above] {$(g,-)$};
    				\draw [->,thick] (X2) -- (X3) node [midway,above] {$(f,-)$};
    				\draw [->,thick] (X2) -- (X4) node [midway,left] {$t$};
    				\draw [->,thick,dashed] (X3) -- (X4) node [midway,right] {$s$};
    				\draw [->,thick] (X1) -- (X4) node [midway,left] {$0$};
    			\end{tikzpicture}
    		\end{center}
    		for any morphism $t$ that makes the left-hand triangle commutative, there exists a morphism $s$ that makes the right-hand triangle commutative.
    		\item[(2)]
    		$\mathbb{X}=X\overset{f}{\rightarrow} Y\overset{g}{\rightarrow}Z$ is not necessarily a complex, if it is indeed a complex, i.e. $gf=0$, then
    		\begin{center}
    			$\Omega_{\mathbb{X}}:=\{M\in \Mod\text{-}\mathcal{X} \mid M(X)\overset{M(f)}{\longrightarrow} M(Y)\overset{M(g)}{\longrightarrow}M(Z)\quad \text{is exact}\}$.
    		\end{center}
    	\end{itemize}
    \end{remark}
    
    \begin{theorem}\label{3.9}
    	A subcategory $\mathscr{D}$ of $\Mod$-$\mathcal{X}$ is definable if and only if it satisfies the following condition.
    	\begin{itemize}
    		\item There is a family $(\mathbb{X}_{\lambda}=X_{\lambda}\overset{f_{\lambda}}{\rightarrow} Y_{\lambda}\overset{g_{\lambda}}{\rightarrow} Z_{\lambda})_{\lambda\in \Lambda}$ of objects from $\Romanbar{3}(\mathcal{X})$ such that 
    		\begin{center}
    			$\mathscr{D}=\displaystyle{\bigcap_{\lambda\in \Lambda}}\Omega_{\mathbb{X}_{\lambda}}$.
    		\end{center}
    	\end{itemize}
       \begin{proof}
       	First let $\mathscr{D}$ be a definable subcategory of $\Mod\text{-}\mathcal{X}$. Therefore, there is a family of coherent functors $(\mathbb{F}_{\lambda})_{\lambda\in \Lambda} \in\modd\text{-}(\modd\text{-}\mathcal{X})$, such that
       	$\mathscr{D}=\displaystyle{\bigcap_{\lambda\in \Lambda}}\Ker(\overrightarrow{\mathbb{F}}_{\lambda})$. If we put $\mathbb{X}_{\lambda}:=\mathfrak{k}(\mathbb{F}_{\lambda})$, by Lemma \ref{3.5} $\mathscr{D}=\displaystyle{\bigcap_{\lambda\in \Lambda}}\Omega_{\mathbb{X}_{\lambda}}$.
       	
       	Conversely, let $\mathscr{D}=\displaystyle{\bigcap_{\lambda\in \Lambda}}\Omega_{\mathbb{X}_{\lambda}}$ for a family $(\mathbb{X}_{\lambda}=X_{\lambda}\overset{f_{\lambda}}{\rightarrow} Y_{\lambda}\overset{g_{\lambda}}{\rightarrow} Z_{\lambda})_{\lambda\in \Lambda}$ of objects from $\Romanbar{3}(\mathcal{X})$. If we put $\mathbb{F}_{\lambda}:=\mathfrak{k}^{-1}(\mathbb{X}_{\lambda})$, then by Lemma \ref{3.6}
       	\begin{center}
       		$\mathscr{D}=\displaystyle{\bigcap_{\lambda\in \Lambda}}\Omega_{\mathbb{X}_{\lambda}}=\displaystyle{\bigcap_{\lambda\in \Lambda}}\Ker(\overrightarrow{\mathbb{F}}_{\lambda})$.
       	\end{center}
       \end{proof}
    \end{theorem}

    In the following (well known) example, we see that most famous definable subcategories are included in the characterization of Theorem \ref{3.9}, and actually Theorem \ref{3.9} is saying that general definable subcategories are not too far from the example.
    
    \begin{example}\label{3.10}
    	Let $\mathcal{X}$ be a skeletally small exact category. Then the following subcategories of $\Mod\text{-}\mathcal{X}$ are definable.
    	\begin{itemize}
    		\item[(1)] The category $\rm Lex\text{-}\mathcal{X}$ of left exact functors in $\Mod\text{-}\mathcal{X}$.
    		\item[(2)] The category $\rm Rex\text{-}\mathcal{X}$ of right exact functors in $\Mod\text{-}\mathcal{X}$.
    		\item[(3)] The category $\rm Ex\text{-}\mathcal{X}$ of exact functors in $\Mod\text{-}\mathcal{X}$.
    	\end{itemize}
    \end{example}

    In the end of this section we want to apply the above results to modules over rings. Let $A$ be a ring. First we fix the following notations.
    \begin{itemize}
    	\item[(1)]
    	$\Mod\text{-}A$ is the category of all left $A$-modules.
    	\item[(2)]
    	$\Mod\text{-}A^{op}$ is the category of all right $A$-modules.
    	\item[(3)]
    	$\modd\text{-}A$ is the category of all finitely presented left $A$-modules.
    	\item[(4)]
    	$\modd\text{-}A^{op}$ is the category of all finitely presented right $A$-modules.
    	\item[(5)]
    	$\proj(A)$ is the category of all finitely generated projective left $A$-modules.
    	\item[(4)]
    	$\proj(A^{op})$ is the category of all finitely generated projective right $A$-modules.
    \end{itemize}
    
     Since $\Mod\text{-}A\cong\Mod\text{-}(\proj(A))$ and  $\modd\text{-}A\cong\modd\text{-}(\proj(A))$ we can apply our results to modules over $A$. 
     
     The following notation is the analogous of Notation \ref{3.7} for modules over rings.
     
     \begin{notation}\label{3.11}
     	Let $\mathbb{P}=P\overset{f}{\rightarrow} Q\overset{g}{\rightarrow}R$ be a diagram in additive category $\proj(A)$. Define
     	\begin{center}
     		$\Omega_{\mathbb{P}}:=\{M\in \Mod\text{-}A \mid \text{any morphism}\; t:Q\rightarrow M\; \text{with}\, tf=0\; \text{factors through}\, g\}$.
     	\end{center}
     \end{notation}
 
     	It is clear that if we consider $A$ as a pre-additive category with a single object, the composition of obvious functors $A\rightarrow \proj(A^{op})\rightarrow \Romanbar{3}(\proj(A^{op}))$ corresponds to Adelman's free abelian category. The other presentation is given by
     	\begin{align*}
     		A&\longrightarrow \big[\modd\text{-}(\modd\text{-}A)\big]^{op}.\\
     		A&\longmapsto (A,-)
     	\end{align*}
     The appearance of "$^{op}$", which may seems confusing at first glance, comes from the fact that the embedding functor from $A$, considered as a category with a single object, to the category of left modules is contravariant! 
     
     Using Theorem \ref{3.9} we obtain the following characterization of definable subcategories of $\Mod\text{-}A$.
     
     \begin{theorem}\label{3.12}
     	Let $\mathscr{D}$ be a subcategory of $\Mod\text{-}A$. The following are equivalent.
     	\begin{itemize}
     		\item[(1)]
     		$\mathscr{D}$ is a definable subcategory.
     		\item[(2)]
     		 There is a family $(\mathbb{P}_{\lambda}=P_{\lambda}\overset{f_{\lambda}}{\rightarrow} Q_{\lambda}\overset{g_{\lambda}}{\rightarrow} R_{\lambda})_{\lambda\in \Lambda}$ of objects from $\Romanbar{3}(\proj(A))$ such that
     		 \begin{center}
     		 	$\mathscr{D}=\displaystyle{\bigcap_{\lambda\in \Lambda}}\Omega_{\mathbb{P}_{\lambda}}$.
     		 \end{center}
     		\item[(3)] There is a family of pairs of matrices $(U_{\lambda},V_{\lambda})_{\lambda\in \Lambda}$ (of appropriate size) with entries from $A$ such that 
     		\begin{align}
     			\mathscr{D}=\{M\in \Mod\text{-}A \mid \forall \lambda\in \Lambda \;\text{and}\; \forall x\in M^n, U_{\lambda}x=0\implies x=V_{\lambda}y\; \text{for some}\; y\in M^m\}.\label{MT}
     		\end{align} 
     	\end{itemize}
       \begin{proof}
       	  The equivalence of $(1)$ and $(2)$ follows from Theorem \ref{3.9}.
       	  For the equivalence of $(2)$ and $(3)$, it is enough to replace projective modules, in projective presentations, with free modules of finite rank. Just note that the corresponding pair of matrices for the object
       	 \begin{align*}
       	 	\mathbb{F}= A^m\overset{U}{\rightarrow} A^n\overset{V}{\rightarrow} A^k
       	 \end{align*}
         is $(U^T,V^T)$.
       \end{proof}
     \end{theorem}
      
      For any definable subcategory $\mathscr{D}$ of $\Mod\text{-}A$, there is a family of pairs of matrices $(U_{\lambda},V_{\lambda})_{\lambda\in \Lambda}$ such that the equality \eqref{MT} is satisfied. In this case we say that $\mathscr{D}$ {\it is defined by the family $(U_{\lambda},V_{\lambda})_{\lambda\in \Lambda}$}.
      
      In the following example, we illustrate our characterization, Theorem \ref{3.12}, for a simple definable subcategory.
      
      \begin{example}
      	Let $\mathbb{Z}$ be the ring of integers. Then we can check that
      	\begin{align*}
      		\mathscr{D}:=\Add(\mathbb{Z}/{2\mathbb{Z}})=\{M\in \Mod\text{-}\mathbb{Z} \mid\; 2.x=0,\; \forall x\in M\}
      	\end{align*}
       is a definable subcategory of $\Mod\text{-}\mathbb{Z}$. Indeed if we let $\mathbb{F}:\modd\text{-}\mathcal{X}\rightarrow \Ab$ be the finitely presented functor given by projective presentation
       \begin{equation*}
       	\Hom(\mathbb{Z}/{2\mathbb{Z}},-)\rightarrow\Hom(\mathbb{Z},-)\rightarrow \mathbb{F}\rightarrow 0,
       \end{equation*} 
       then $\mathscr{D}=\Ker(\overrightarrow{\mathbb{F}})$. If we apply Theorem \ref{3.12}, we have that $\mathscr{D}=\Omega_{\mathbb{P}}$ where $\mathbb{P}$ is the following object of $\Romanbar{3}(\proj(\mathbb{Z}))$.
       \begin{align*}
       	\mathbb{Z}\overset{[-1,\;2]^T}{\longrightarrow} \mathbb{Z}\oplus\mathbb{Z}\overset{[0,\;-1]}{\longrightarrow} \mathbb{Z}
       \end{align*}
      and so $\mathscr{D}$ is defined by the pair of matrices $([-1,\;2],[0,\;-1]^T)$.
      \end{example}
    
    \section{Auslander-Gruson-Jensen duality}
    In this section, using the results of the Section 2, we give a conceptual proof of Auslander-Gruson-Jensen duality \cite{GJ,Aus86}. And show how this interpretation makes the duality between definable subcategories of left and right modules, due to Herzog \cite{Her93} for modules over rings, more transparent. First let recall Auslander-Gruson-Jensen duality in a general form as in \cite[Section 10.3]{Pre09}.
    
    Let $\mathcal{X}$ be an additive category. By Theorem \ref{2.5}
    \begin{align*}
    	\mathfrak{i}:&\mathcal{X}\longrightarrow \modd\text{-}(\modd\text{-}\mathcal{X})\\
    	&X\longmapsto \big((X,-),-\big)
    \end{align*}
    and
    \begin{align*}
    	\mathfrak{i}^{op}:&\mathcal{X}^{op}\longrightarrow \modd\text{-}(\modd\text{-}\mathcal{X}^{op})\\
    	&X\longmapsto \big((-,X),-\big)
    \end{align*}
    provide us free abelian categories over $\mathcal{X}$ and $\mathcal{X}^{op}$, respectively. By the universal property of the free abelian category we have:
    
    \begin{theorem}\label{4.1}
    	There is an exact duality $d:\modd\text{-}(\modd\text{-}\mathcal{X}) \longrightarrow \modd\text{-}(\modd\text{-}\mathcal{X}^{op})$ which makes the following diagram commutative, where $\mathfrak{d}$ is the trivial duality between $\mathcal{X}$ and $\mathcal{X}^{op}$.
    	\begin{center}
    		\begin{tikzpicture}
    			\node (X0) at (-1,1.5) {$\mathcal{X}$};
    			\node (X1) at (-1,-1.5) {$\mathcal{X}^{op}$};
    			\node (X2) at (3,1.5) {$\modd\text{-}(\modd\text{-}\mathcal{X})$};
    			\node (X3) at (3,-1.5) {$\modd\text{-}(\modd\text{-}\mathcal{X}^{op})$};
    			\draw [->,thick] (X0) -- (X1) node [midway,left] {$\mathfrak{d}$};
    			\draw [->,thick] (X0) -- (X2) node [midway,above] {$\mathfrak{i}$};
    			\draw [->,thick] (X1) -- (X3) node [midway,above] {$\mathfrak{i}^{op}$};
    			\draw [->,thick,dashed] (X2) -- (X3) node [midway,right] {$d$};
    		\end{tikzpicture}
    	\end{center}
    \end{theorem}
    \begin{remark}\label{4.2}
    	The exact duality $d:\modd\text{-}(\modd\text{-}\mathcal{X}) \longrightarrow \modd\text{-}(\modd\text{-}\mathcal{X}^{op})$ is known as Auslander-Gruson-Jensen duality. Its definition is somewhat complicated. But using \ref{dia2.2} we can easily define it as the unique exact contravariant functor that sends $\big((X,-),-\big)\rightarrow \big((Y,-),-\big)$ to $\big((-,X),-\big)\leftarrow \big((-,Y),-\big)$. 
    \end{remark}
Also by Theorem \ref{2.3}
\begin{align*}
	\mathfrak{j}:&\mathcal{X}\longrightarrow \Romanbar{3}(\mathcal{X})\\
	&X\longmapsto \big(0\rightarrow X\rightarrow 0)
\end{align*}
and
\begin{align*}
	\mathfrak{j^{op}}:&\mathcal{X}^{op}\longrightarrow \Romanbar{3}(\mathcal{X}^{op})\\
	&X\longmapsto \big(0\leftarrow X\leftarrow 0)
\end{align*}
provide us free abelian categories over $\mathcal{X}$ and $\mathcal{X}^{op}$, respectively. Denoting by $d_{\scriptsize\Romanbar{3}}:\Romanbar{3}(\mathcal{X})\rightarrow \Romanbar{3}(\mathcal{X}^{op})$ the trivial duality obtained by reversing all arrows, we have:

\begin{theorem}\label{4.3}
	$d_{\scriptsize\Romanbar{3}}$ is the unique exact contravariant functor that makes the following diagram commutative.
	\begin{center}
		\begin{tikzpicture}
			\node (X0) at (-1,1.5) {$\mathcal{X}$};
			\node (X1) at (-1,-1.5) {$\mathcal{X}^{op}$};
			\node (X2) at (3,1.5) {$\Romanbar{3}(\mathcal{X})$};
			\node (X3) at (3,-1.5) {$\Romanbar{3}(\mathcal{X}^{op})$};
			\draw [->,thick] (X0) -- (X1) node [midway,left] {$\mathfrak{d}$};
			\draw [->,thick] (X0) -- (X2) node [midway,above] {$\mathfrak{j}$};
			\draw [->,thick] (X1) -- (X3) node [midway,above] {$\mathfrak{j}^{op}$};
			\draw [->,thick,dashed] (X2) -- (X3) node [midway,right] {$d_{\scriptsize\Romanbar{3}}$};
		\end{tikzpicture}
	\end{center}
\end{theorem}

\begin{remark}\label{4.4}
	By the uniqueness of free abelian category, we can see $d_{\scriptsize\Romanbar{3}}$ as another interpretation of Auslander-Gruson-Jensen duality. Indeed we have the following commutative diagram of functors.
	\begin{center}
		\begin{tikzpicture}
			\node (X1) at (-4,1) {$\Romanbar{3}(\mathcal{X})$};
			\node (X2) at (0,1) {$\Romanbar{3}(\mathcal{X}^{op})$};
			\node (X3) at (-4,-1) {$\modd$-$(\modd$-$\mathcal{X})$};
			\node (X4) at (0,-1) {$\modd$-$(\modd$-$\mathcal{X}^{op})$};
			\draw [->,thick] (X1) -- (X2) node [midway,above] {$d_{\scriptsize\Romanbar{3}}$};
			\draw [->,thick] (X1) to [out=240,in=120] (X3) node [midway,above] {};
			\draw [->,thick] (X2) to [out=240,in=120] (X4) node [midway,left] {};
			\draw [->,thick] (X3) to [out=60,in=300] (X1) node [midway,above] {};
			\draw [->,thick] (X4) to [out=60,in=300] (X2) node [midway,above] {};
			\draw [->,thick] (X3) -- (X4) node [midway,above] {$d$};
			\node (K1) at (-4.6,0) {$\mathfrak{k}$};
			\node (K2) at (-3.2,0) {$\mathfrak{k}^{-1}$};
			\node (K3) at (-0.6,0) {$\mathfrak{k}$};
			\node (K4) at (0.7,0) {$\mathfrak{k}^{-1}$};
		\end{tikzpicture}
	\end{center}
\end{remark}

    Herzog \cite{Her93} proved that for a ring $A$, there is a natural bijection between definable subcategories of $\Mod\text{-}A$ and definable subcategories of $\Mod\text{-}A^{op}$ (see also \cite{Baz08}). Using our characterization of definable subcategories and the above presentation of Auslander-Gruson-Jensen duality, we prove a general form of this bijection, in a more illuminating way.
    
    Let $\mathbb{X}=X\overset{f}{\rightarrow} Y\overset{g}{\rightarrow}Z$ be an object in $\Romanbar{3}(\mathcal{X})$. We defined in Notation \ref{3.7}
    \begin{center}
    	$\Omega_{\mathbb{X}}=\{M\in \Mod\text{-}\mathcal{X} \mid \Ker(M(g))\subseteq \Imm(M(f))\}$.
    \end{center}
    Similarly, we define
    \begin{center}
    	$\Gamma_{\mathbb{X}}=\{N\in \Mod\text{-}\mathcal{X}^{op} \mid \Ker(N(f))\subseteq \Imm(N(g))\}$.
    \end{center}
    By Theorem \ref{3.9} every definable subcategory of $\Mod\text{-}\mathcal{X}$ is of the form
    \begin{center}
    	$\mathscr{D}=\displaystyle{\bigcap_{\lambda\in \Lambda}}\Omega_{\mathbb{X}_{\lambda}}$,
    \end{center}
    for a family $(\mathbb{X}_{\lambda})_{\lambda\in \Lambda}$ of objects in $\Romanbar{3}(\mathcal{X})$.
    Then if we put
    \begin{center}
    	$\bar{\mathscr{D}}=\displaystyle{\bigcap_{\lambda\in \Lambda}}\Gamma_{\mathbb{X}_{\lambda}}$,
    \end{center}
    clearly we obtain a definable subcategory of $\Mod\text{-}\mathcal{X}^{op}$. More precisely
    
    \begin{theorem}\label{4.5}
    	Sending a definable subcategory $\mathscr{D}\subseteq\Mod\text{-}\mathcal{X}$ to $\bar{\mathscr{D}}$, gives a bijection between definable subcategories of $\Mod\text{-}\mathcal{X}$ and definable subcategories of $\Mod\text{-}\mathcal{X}^{op}$.
    \end{theorem}

In the following we apply Theorem \ref{4.5} to definable subcategories in Example \ref{3.10}. 

    \begin{example}
    		Let $\mathcal{X}$ be a skeletally small exact category. Then the bijection of Theorem \ref{4.5} sends 
    	\begin{itemize}
    		\item[(1)]
    		The subcategory $\rm Lex\text{-}\mathcal{X}\subseteq \Mod\text{-}\mathcal{X}$ of left exact functors to $\rm Rex\text{-}\mathcal{X}^{op}\subseteq \Mod\text{-}\mathcal{X}^{op}$, the subcategory of contravariant right exact functors.
    		\item[(2)]
    		The subcategory $\rm Rex\text{-}\mathcal{X}\subseteq \Mod\text{-}\mathcal{X}$ of right exact functors to $\rm Lex\text{-}\mathcal{X}^{op}\subseteq \Mod\text{-}\mathcal{X}^{op}$, the subcategory of contravariant left exact functors.
    		\item[(3)]
    		The category $\rm Ex\text{-}\mathcal{X}\subseteq \Mod\text{-}\mathcal{X}$ of exact functors to $\rm Ex\text{-}\mathcal{X}^{op}\subseteq \Mod\text{-}\mathcal{X}^{op}$, the subcategory of contravariant exact functors.
    	\end{itemize}
    \end{example}
For the sake of completeness, we state the above results for the case of modules over a ring $A$. The proofs are (again) left to the reader.
\begin{corollary}
	The duality $()^*=\Hom(-,A):\Romanbar{3}(\proj(A))\rightarrow \Romanbar{3}(\proj(A^{op}))$, makes the following diagram commutative, where the bottom row is the Auslander-Gruson-Jensen duality and vertical functors are canonical exact equivalences.
	\begin{center}
		\begin{tikzpicture}
			\node (X1) at (-4,1) {$\Romanbar{3}(\proj(A))$};
			\node (X2) at (0,1) {$\Romanbar{3}(\proj(A)^{op})$};
			\node (X3) at (-4,-1) {$\modd$-$(\modd$-$A)$};
			\node (X4) at (0,-1) {$\modd$-$(\modd$-$A^{op})$};
			\draw [->,thick] (X1) -- (X2) node [midway,above] {};
			\draw [->,thick] (X1) -- (X3) node [midway,above] {};
			\draw [->,thick] (X2) -- (X4) node [midway,left] {};
			\draw [->,thick] (X3) -- (X4) node [midway,left] {};
		\end{tikzpicture}
	\end{center}
\end{corollary}

   Let $\mathscr{D}$ be a definable subcategory of $\Mod\text{-}A$. By Theorem \ref{3.12} there is a family of pairs of matrices $(U_{\lambda},V_{\lambda})_{\lambda\in \Lambda}$ (of appropriate size) with entries from $A$ such that 
   \begin{align*}
   	\mathscr{D}=\{M\in \Mod\text{-}A \mid \forall \lambda\in \Lambda \;\text{and}\; \forall x\in M^n, U_{\lambda}x=0\implies x=V_{\lambda}y\; \text{for some}\; y\in M^m\}.
   \end{align*} 
    And in this case we said that $\mathscr{D}$ is defined by the family $(U_{\lambda},V_{\lambda})_{\lambda\in \Lambda}$. Clearly the family $(V_{\lambda}^T,U_{\lambda}^T)_{\lambda\in \Lambda}$ defines a definable subcategory in $\Mod\text{-}A^{op}$. Indeed we have
       \begin{corollary}
             Sending a definable subcategory $\mathscr{D}$ of $\Mod\text{-}A$ defined by a family $(U_{\lambda},V_{\lambda})_{\lambda\in \Lambda}$ of pairs of matrices with entries in $A$, to the definable subcategory $\bar{\mathscr{D}}$ of $\Mod\text{-}A^{op}$ defined by the family $(V_{\lambda}^T,U_{\lambda}^T)_{\lambda\in \Lambda}$, gives a bijection between definable subcategories of $\Mod\text{-}A$ and definable subcategories of $\Mod\text{-}A^{op}$.
       \end{corollary}

	\section*{acknowledgements}
	The author appreciates Professor Mike Prest for his valuable feedback and suggestions.
	The research was supported by a grant from IPM.


\end{document}